\documentclass[12pt]{article}
\usepackage{amssymb}
\usepackage{makeidx}
\usepackage{amsbsy,amsmath}
\usepackage{latexsym,amssymb,mathptm}
\usepackage{bm}

\usepackage[usenames]{color}

\newtheorem{Theorem}{Theorem}

\title{{\bf The Distribution of Argmaximum or \\a Winner Problem} }
\author{{Youri \sc Davydov${}^{1}$}\quad  and\quad { Vladimir \sc Rotar ${}^{2}$}}
\date
{\footnotesize  ${}^{1}$  Laboratoire Paul Painlev\'e \\
 Universit\'e de Lille, France,\\
 and Faculty of Mathematics and Computer Sciences\\
 St. Petersburg State University, Russia\\
 Email: youri.davydov@univ-lille.fr\\
\vspace{3pt}
${}^{2}$ \ Department of Mathematics \\
 University of California at San Diego, USA and \\
 the National University, USA \\
  Email: vrotar@math.ucsd.edu }

\begin{document}

\maketitle

%\begin{abstract}
\begin{quote}
{\bf Abstract.}
We consider a limit theorem for the distribution of a random variable (r.v.) $A_n=\arg\max_{i: 1\dots n}\{X_i\}$, where $X_i$'s are independent continuous non-negative random variables. The r.v.'s $X_i,\; i=1.\dots, n$, may be interpreted as the gains of $n$ players in a game, and the   r.v. $Y_n$ itself as the number of a ``winner". In the case of independent identically distributed r.v.'s, the distribution of $A_n$  is, clearly, uniform on $\{1,\dots,  n\}$, while when the $X$'s are non-identically distributed, the problem is not trivial. The paper contains some limit theorems for the distribution of $A_n$ as $n\to\infty$.
\end{quote}


\bigskip \bigskip \noindent Keywords: limit theorem, maximum of random variables, argmaximum.
%\end{abstract}
%\end{quote}

\section{Introduction}

\renewcommand{\theequation}{\thesection.\arabic{equation}}
\setcounter{equation}{0}
%\subsection{Background and Motivation}

\setcounter{equation}{0}

Let  $X_{1},X_{2},...$ be  positive  and  independent random variables (r.v.'s). We will deal with $\max\{X_1,...,X_n\}$.  For the case of identically distributed r.v.'s, the  theory of limiting distribution for the maximum was developed in  papers by Fisher\&Tippett \cite{fisher}, von Mises \cite{mises} and Gnedenko \cite{gnedenko}; see also systematic presentations  in \cite{feller},  \cite{haan},  \cite{rotar1}, \cite{rotar2}.

The case of non-identically distributed r.v.'s was being studied not so thoroughly. We mention here three works: \cite{rotar1} that concerns the merging of  the distributions of the maximum for two triangular series of the $X$'s, and  \cite{eks} and  \cite{s}  that are devoted to the marginal maxima of Gaussian vectors.

This paper concerns the distribution of r.v. $A_n= \arg\max\{X_1,\ldots,X_n\}.$
{If the r.v.'s} $X_i,\; i=1.\dots, n$, { are interpreted as the gains of} $n$ {players in a game, then} $A_n$ {is the number of  a ``winner".}

  Below, we assume the  $X$'s to be continuous, and in this case, the definition of $ \arg\max\{X_1,\ldots,X_n\}$ is correct { (the winner is almost surely unique)}.

  We believe that, having its intrinsic value, the results concerning limiting behavior of   $A_n$ may be helpful in applications too.

   As an example, consider  a complex machine consisting of a large number of parts with random and non-identically distributed lifetimes. The question is which part will break down first. Certainly, we deal here with $\arg\min $ but it can be easily reduced to $\arg\max$. \footnote{We thank  professor Vadim Ponomarenko (San Diego State University) for this example.}

Let
\[
p_{in}: = P\{A_n=i\}= P\{X_{i}= \max\{X_1,...,X_n\}\}, \,\,\,\,  i=1,\dots ,n.
\]

In the case of i.i.d. r.v.'s, the probability   $p_{in}$ is, clearly, equal to $1/n$,  while {the case where  $X$'s are non-identically distributed, requires -- as we will see --  some work.}

Let  $F_{i}(x)=P(X_{i}\leq x),  \,\,\, F(0)=0, F(x)>0$ for $x>0$,
  Set for $x>0$
\[
\nu_i(x)=-\ln F_i(x). \,\,\,\,
\]
and $\nu_i(0)=\infty$.

 So, for all $i$,
 \begin{eqnarray}
 && F_i(x)=\exp\{-\nu_i(x)\},\label{nu0}   \\
 && \nu_i(x) \,\,\,\text{ is non-increasing},     \,\,\, \nu_i(0)=\infty, \,\,\, \nu_i(\infty)=0.  \label{nu1}
 \end{eqnarray}

 The asymptotic behavior of $\nu_i(x) $ as $x\to\infty$ is equivalent to that of $1-F_i(x)$.\vspace{.1in}
Since we assumed   all distributions $F_i$ to be  non-atomic,  $\nu_i(x)$'s are continuous for $x>0$.

{As a \textit{basic example}, one may consider a family of distributions that satisfies the following condition:}

\begin{equation}\label{r1}
\nu_i(x)=c_ir(x),
\end{equation}
where $r(x)$ is a non-negative, continuous,  and non-increasing function; $r(0)=\infty$, \, $r(\infty)=0$, and $c_i$'s are non-negative.

\vspace{.1in}

%\textit{\textbf{Example} (\textit{Weibull's distrinution})   }

Well known Weibull's distrinutions
\[
   F_i(x)=\exp\left\{-\frac{c_i}{x^\alpha}\right\},\;\; \text{and} \,\,\, c_i, \,\,\alpha>0,
\]
which are stable with respect to maximization, may serve as a particular case of such a family.

%\textit{\textbf{Basic Example} \[\{bfrepresent an important example of a family  which satisfies the following condition:

%$\square$ \footnote{the symbol $\square$  means the end of an example; the symbol $\blacksquare$ below will mean  the end of a proof. } \vspace{.1in}

We will show below that under condition (\ref{r1}),

\[
p_{in}= \frac{c_i}{\sum_{i=1}^nc_i}. \,\,\,
\]

However, our main  goal is to prove  that for large $n$, to get a similar representation, we may proceed just from the asymptotic behavior of $\nu_i(x)$ (or $1-F_i(x)$) at infinity.

\section{Main Results }\label{subsecprop}

\renewcommand{\theequation}{\thesection.\arabic{equation}}
\setcounter{equation}{0}

Assume the following.
\begin{enumerate}
\item
 \begin{equation}\label{cond1}
 \nu_i(x)=c_i r(x) (1+\delta_i(x)),
  \end{equation}
where $r(x)$ is defined as above, $\delta_i(x)$ are continuous,  uniformly in $i$
 \begin{equation}\label{cond1+}
   \delta_i(x)\to 0 \,\,\, \text{as}\,\,\, x\to\infty,
 \end{equation}
and for  positive constants $M<\infty$ and $m<1$, and for all $i$ and $x$,
\begin{equation}\label{cond1++}
-m\leq \delta_i(x)\leq M.
\end{equation}
  \item
    \begin{equation}\label{cond2}
     b_n=:\sum_{i=1}^nc_i \to\infty \,\,\, \text{ as }\,\,\,  n\to \infty.
     \end{equation}
      \end{enumerate}

\subsection{A local limit theorem}\label{localsection}

\begin{Theorem}
\label{th1}
 Set
\[
\alpha_{in}=\frac{c_{i}}{b_n}.
\]
Then, under the above conditions, uniformly in  $i$
\begin{equation} \label{p1}
p_{in}\sim \alpha_{in} \,\,\,\;\;\text{as}\,\,\,n\to\infty.
\footnote{The symbol $\sim$ means that the ratio of the left- and right-hand sides converges to one.}
\end{equation}
\end{Theorem}
\vspace{.1in}
\subsection{An integral  limit theorem}\label{lintegralsection}

To obtain an integral limit theorem, we should suppose that the coefficients $c_i$'s  are varying --  in a sense -- regularly.

Consider the segment $[0,1]$ and identify a point $i/n, \,\,\,i=1.\dots,n$,   with a r.v. $X_i$; so to speak,  with the $i$-th ``player". (In other words, we consider r.v.\ $A_n/n$.)

Let us assign to the point $i/n$    probability  $\alpha_{in}$, and suppose that the measure so defined weakly  converges to a probability measure $\alpha$ on $[0,1]$. In other terms,
\begin{equation}\label{alpha}
	\alpha_n=: \sum_{i=1}^n\delta_{\{i/n\}}\alpha_{in}\Rightarrow\alpha,
\end{equation}
where $\delta_{\{x\}}$ is a measure concentrated at point $x$. We will discuss condition (\ref{alpha}) in detail in Section \ref{rv}.  \vspace{0.1in}

\vspace{.1in}
\begin{Theorem}
	\label{th2}
Suppose that together with conditions of Theorem \ref{th1}, (\ref{alpha}) holds. Then discrete measure
\begin{equation}\label{alpha1}
 \mu_n=: \sum_{i=1}^n\delta_{\{i/n\}}p_{in}\Rightarrow\alpha,
\end{equation}
 \end{Theorem}
\vspace{.1in}
{\bf Proof.} Theorem \ref{th2} easily follows from Theorem \ref{th1}. Indeed, since the convergence in (\ref{p1})  is uniform in $i$, for any continuous bounded  function $h$,
\begin{eqnarray*}
	% \nonumber to remove numbering (before each equation)
	\int_0^1hd\mu_n &=&\sum_1^nh\left(\frac{i}{n}\right)p_{in}=\sum_1^nh\left(\frac{i}{n}\right)\alpha_{in}(1+o(1)) \\
	&=& (1+o(1))\int_0^1hd\alpha_n \to\int_0^1hd\alpha. \,\,\,\, \,\,\,\,
\end{eqnarray*}
\begin{flushright}
	$\blacksquare$
\end{flushright}
\vspace{.1in}

\section{Proofs}

\renewcommand{\theequation}{\thesection.\arabic{equation}}
\setcounter{equation}{0}

\subsection{A basic formula}\label{basicformula}
 We have
 \begin{eqnarray*}
  p_{in}&=& \int_{0}^{\infty }\prod_{j=1,\, j\neq i}^{n}\,F_{j}(x)dF_{i}(x)\\ \nonumber
  &=& - \int_0^\infty \exp\left\{-\sum_{j=1, j\neq i}^n \nu_j(x) \right\}\exp\{-\nu_i(x)\}d\nu_i(x) \nonumber \\
  &=& -  \int_0^\infty \exp\left\{-\sum_{j=1}^n \nu_j(x) \right\}d \nu_i(x).
\end{eqnarray*}

Integrating by parts  and taking into account (\ref{nu1}), we get
\begin{equation}
  p_{in}=    -  \int_0^\infty \nu_i(x)  \exp\left\{-\sum_{j=1}^n \nu_j(x) \right\}d\left(\sum_{j=1}^n \nu_j(x)\right). \label{nu2}
\end{equation}

Consider substitution
\begin{equation}\label{sub}
\sum_{i=1}^n \nu_i(x)=y.
\end{equation}
For any non-increasing function $r(x)$, we define its inverse as
\[
r^{-1}(y)=\sup\{x: r(x)\geq y\}.
\]
Let $x_n(y)$ be the inverse of the function $\sum_{i=1}^n \nu_i(x)$; in other words, a solution (in the above sense) to  equation (\ref{sub}).
Then from (\ref{nu2})--(\ref{sub}), it follows that

\begin{equation}\label{p-in}
    p_{in}= \int_0^\infty \nu_i(x_n(y)) e^{-y}dy.
\end{equation}

This may serve as a basic formula. \vspace{.1in}

\textbf{Remark. } Condition  $F_i(x)>0$ for all $x>0$  is not necessary; we imposed it just to make the proof of (\ref{p-in}) more explicit.  As a matter of fact, it is easy (though a bit longer) to prove that the same is true, for example, if  for all $n$ and a finite  $a\geq0$
 \[
 a_n=:  \max_{i=1,\dots n}\sup\{x: F_i(x)=0\}\leq a.
 \]

To better understand how   formula (\ref{p-in}) may work, consider first the basic example   (\ref{r1}) where  
  $r(x)$ is a non-negative, continuous,  and non-increasing function; $r(0)=\infty$, \, $r(\infty)=0$, and $c_i$'s are non-negative.
Then
\[
\sum_{i=1}^n \nu_i(x)=r(x)\sum_{i=1}^n c_i,
\]
and a solution to equation (\ref{sub}) is
 \begin{equation}\label{x_n-example}
  x_n(y)=r^{-1}\left(\frac{y}{\sum_{i=1}^nc_i}\right).
 \end{equation}
So,
\begin{equation}\label{r2}
\nu_i(x_n(y))=c_ir\left(r^{-1}\left(\frac{y}{\sum_{i=1}^nc_i}\right)\right) =\frac{c_i}{\sum_{i=1}^nc_i}\,y.
\end{equation}
Thus, in this case,
\[
p_{in}=\frac{c_i}{\sum_{i=1}^nc_i}\int_0^\infty ye^{-y}dy=\frac{c_i}{\sum_{i=1}^nc_i}.
\]

When looking at (\ref{x_n-example}), one may suppose that for large $n$,  the asymptotic behavior of $p_{in}$ is based just on the asymptotics  of $r^{-1}(x)$ at zero, which is connected with that of $r(x)$ at infinity (or tails $1-F_i(x))$.
\subsection{Proof of Theorem \ref{th1}}

     Let $x_n(y)$ be  a solution to equation
     \begin{equation}\label{eq1}
        \sum_{i=1}^n \nu_i(x)=y,
     \end{equation}
     that is,
      \begin{equation}\label{eq1+}
        r(x)\sum_{i=1}^n c_i(1+\delta_i(x))=y.
     \end{equation}
So,
 \begin{equation}\label{x_n}
      x_n(y) =r^{-1}\left(\frac{y}{\sum_{i=1}^n c_i(1+\delta_i(x_n(y)))}\right) .
 \end{equation}

From (\ref{cond1++}),  it follows that

\begin{equation}\label{x_n1}
x_n(y) \geq  r^{-1}\left(\frac{y}{(1-m)\sum_{i=1}^n c_i}\right).
\end{equation}
Hence, since $r^{-1}(0)=\infty$, and  in view of (\ref{cond2}),
\begin{equation}\label{xtoinfty}
x_n(y)\to \infty
\end{equation}
as $n\to\infty$.

 Furthermore, in view of (\ref{x_n}),
     \[
     \nu_i(x_n(y))=c_i\cdot\frac{y(1+\delta_i(x_n(y)))}{\sum_{j=1}^nc_j(1+\delta_j(x_n(y)))}.
     \]
  Thus,
\begin{eqnarray*}
                                 p_{in}&=&\int_0^\infty \nu_i(x_n(y))e^{-y}dy \\
                                 &=& \frac{c_i}{\sum_{j=1}^n c_j} \int_0^\infty y\cdot \frac{(1+\delta_i(x_n(y)))\sum_{j=1}^n c_j}{\sum_{j=1}^n c_j(1+\delta_j(x_n(y)))}\cdot e^{-y}dy
\end{eqnarray*}

   For each $y>0$, in view of (\ref{xtoinfty}) and (\ref{cond1+}),
    \[
    \frac{(1+\delta_i(x_n(y)))\sum_{j=1}^n c_j}{\sum_{j=1}^nc_j(1+\delta_j(x_n(y))}\to 1 \,\,\,\text{as}\,\, n\to\infty,
    \]
   uniformly in $i$. On the other hand,
     \[
    \frac{(1+\delta_i(x_n(y)))\sum_{j=1}^n c_j}{\sum_{j=1}^nc_j(1+\delta_j(x_n(y)))}\leq \frac{1+M}{1-m}.
    \]
   Hence,
   \[
   \frac{p_{in}}{\alpha_{in}}\to \int_0^\infty ye^{-y}dy.\,\,\,\, \,\,\,\,
   \]
\begin{flushright}
	$\blacksquare$
\end{flushright}

%\pagebreak

\section{On regularity of coefficients $c_i$\,\,\,(condition (\ref{alpha})) }\label{rv}
\setcounter{equation}{0}
 In this section we will consider examples illustrating (\ref{alpha}) and present additional results on a possible structure of the limiting distribution $\alpha$ in (\ref{alpha}). 
%\vspace{7pt}
\newpage

\textbf{EXAMPLES.}
\vspace{.1in}
\begin{enumerate}
  \item  Set $c_i=i^s$, for $s\geq0$, and let  $x\in(0,1]$. Let $k=k_n$ be such that $\frac{k}{n}\leq x<\frac{k+1}{n}$.  Then, as is easy to verify,
  \begin{equation}\label{cond3}
   \frac{\sum_{i=1}^{k_n}c_i}{\sum_{i=1}^nc_i}\to x^{s+1}.
   \end{equation}
  (For $k_n=0$, we set  $\sum_{i=1}^{k_n}=0$.) In other words, if $F(x)$ is the distribution function (d.f.) of $\alpha$, then $F(x)=x^{s+1}$.

   Say, if $c_i=i$, then for large $n$, the distribution of the winner numbers may be well presented by a distribution on $[0,1]$ with d.f. $F(x)=x^2$.
  \item Let $c_i=2^i$, Then in the same notations, for any $x<1$,
  \begin{equation}\label{cond3}
   \frac{\sum_{i=1}^{k_n}c_i}{\sum_{i=1}^nc_i}\to 0,
   \end{equation}
   and  measure $\alpha$  is concentrated at point $1$.
     \item Let $c_i=1/i$, Then, as is easy to verify, for any $x\in (0,1]$,
  \begin{equation}\label{cond3}
   \frac{\sum_{i=1}^{k_n}c_i}{\sum_{i=1}^nc_i}\to 1,
   \end{equation}
   and  measure $\alpha$  is concentrated at point $0$. \,\,\, $\square$
\end{enumerate}

As a matter of fact, the class of possible limiting distributions $\alpha$ is narrow because, as we will see,  in (\ref{alpha}) we deal with regularly varying functions (reg.v.f.'s).\footnote{A positive function $H(x)$ on $[0,\infty)$ is regular varying in the sense of Karamata with an order of $\rho, \,\,\,-\infty<\rho<\infty $,  iff for any $x>0$
\[
 \frac{H(tx)}{H(x)}\to x^\rho \,\,\, \text{as}\,\,\, t\to\infty.
 \]
A function $L(x)$ is called slowly varying if it is regularly varying with $\rho=0$. Any reg.v.f.\ $H(x)=x^\rho L(x)$, where $L(\cdot)$ is slowly varying. A detailed presentation of reg.v.f.'s is given, for example, in Feller, \cite{feller}, Chapter VIII, Section 8.  Some definitions and examples may be also found in \cite{rotar1}, Ch,15.  }
\begin{Theorem}
\label{pr3}
\begin{description}
  \item[(A)] Suppose  (\ref{alpha}) holds, and
\begin{equation}\label{regular}
\frac{b_{n+1}}{b_n}\to 1\,\,\,\text{as}\,\,\,n\to\infty.
\end{equation}
Then   the d.f. of $\alpha$ is
\begin{equation}\label{rho}
 F(x)=x^\rho, \,\,\,\,\, x\in[0,1],
\end{equation}
where  $0\leq\rho\leq\infty$, and $b_n=b(n)$, where $b(t)$ is a non-decreasing reg.v.f.

(In (\ref{rho}), if $\rho=0$, then $F(x)=1$ for all $x\in[0,1]$; if $\rho=\infty$, then $F(x)=0$ for all $x<1$.)
  \item[(B)] Vice versa, let   $b_n=b(n)$, where $b(t)$ is a non-decreasing positive  reg.v.f.   Then (\ref{regular}) holds automatically, and (\ref{alpha}) is true with the d.f.\ $F(x)$ of $\alpha$ defined in (\ref{rho}),
\end{description}
\end{Theorem}
\vspace{.1in}
\textbf{Proof}

\textbf{(A)} Let $F_n(x)$ and $F(x)$ be the d.f.'s of measures $\alpha_n$ and $\alpha$, respectively. Then
\begin{equation}\label{f1}
    F_n(x)\to F(x)
\end{equation}
as $n\to \infty$ for all $x$'s that are continuity points of $F(x)$.

Let $b_0=0$, and for all $t\geq0$ function $b(t)=b_n$ if $t\in[n,n+1)$. We will prove that $ b(t)$ is a reg.v.f.

Let us fix a continuity point  $x$, and let an integer $k=k_n$ be such that $\frac{k}{n}\leq x <\frac{k+1}{n}$. Then from  (\ref{f1}) it follows that
\[
    \frac{b_{k_n}}{b_n}\to F(x)\,\,\,\text{as} \,\,\, n\to\infty.
\]
On the other hand, by definition, $b_{k_n}=b(k_n)=b(nx)$, and hence
\begin{equation}\label{f2}
   \frac{b(nx)}{b(n)}\to F(x)\,\,\,\text{as} \,\,\,  n\to\infty.
\end{equation}
Together with (\ref{regular}), this implies that
\begin{equation}\label{f3}
 \frac{b(tx)}{b(t)}\to F(x)\,\,\,\text{as}\,\,\,  t\to\infty,
\end{equation}
where $t$'s are arbitrary positive numbers.
Indeed,  let $n=n_t$ be such that  $t\in[n, n+1)$. Then
 \begin{equation}\label{f4}
    \frac{b(nx)}{b(n+1)}\leq \frac{b(tx)}{b(t)}\leq \frac{b((n+1)x)}{b(n)}
 \end{equation}
\
Furthermore, if $t\to\infty$, then $n=n_t\to\infty$, and
\[
\frac{b(nx)}{b(n+1)}=\frac{b(n)}{b(n+1)}\cdot\frac{b(nx)}{b(n)}\to F(x)
\]
in view of (\ref{regular}) and (\ref{f2}). Similarly,  the same is  true for the very right fraction in (\ref{f4}).

So, function $b(t)$ is  a  regularly varying function, and the limit in  (\ref{f3}) must be equal to a power function  $x^\rho$; see, for instance,  Lemma 1 from Feller \cite{feller}, VIII, 8.

\textbf{(B)} Let $b_n=b(n)$ where $b(t)$  is a reg.v.f.\ (that may be  different from  the piecewise constant function $b(x)$ defined in part (A) of the proof). Let us fix an $x\in(0,1]$, and let again  an integer $k=k_n$ be such that $\frac{k}{n}\leq x <\frac{k+1}{n}$.

First, since $b(x)$ is non-decreasing,
\begin{equation}\label{f5}
 F_n(x)=\frac{b_{k_n}}{b_n}=\frac{b(k_n)}{b(n)}\leq \frac{b(nx)}{b(n)}\to x^\rho,
\end{equation}
 where $0\leq\rho<\infty$.
On the other hand,
\begin{equation}
 F_n(x)=\frac{b(k_n)}{b(n)}\geq \frac{b(xn-1)}{b(n)} =  \frac{b(xn-1)}{b(xn)}\frac{b(xn)}{b(n)}. \label{f6}
\end{equation}
Let us note that for any non-decreasing reg.v.f.\ $b(x)$
\begin{equation}\label{f8}
 \frac{b(x-1)}{b(x)}\to 1  \,\,\, \text{as} \,\,\, x\to\infty.
\end{equation}
 Indeed, for $s<1$ and sufficiently large $x$'s
 \begin{equation*}
  1 \geq  \frac{b(x-1)}{b(x)}\geq  \frac{b(sx)}{b(x)}\to s^\rho,
  \end{equation*}
and the right-hand side can be made arbitrary close to $1$.
 By virtue of (\ref{f8}), the first factor in (\ref{f6}), converges to $1$, and the whole product converges to $x^\rho$.

Relation (\ref{f8}) also implies (\ref{regular}). \,\,\,
\begin{flushright}
	$\blacksquare$
\end{flushright}
\vspace{.1in}

\vspace{.1in}
\section{Conclusive Remarks   }
 \begin{enumerate}
 \item When considering examples, it is more convenient to deal directly with sequences $b_n$ rather than coefficients $c_i$'s. In particular, if $b_n$ are asymptotically   exponential, (\ref{regular}) is not true but it is easy to show that the limiting measure $\alpha $ exists and concentrated at point $1$ (see also Example 2 above). On the other hand, if for instance, $b_n\sim e^{c\sqrt{n}}$ for a positive $c$,   (\ref{regular}) is true though the  limiting measure is again concentrated at $1$.
 \item To specify a particular $\rho$, we may, for example, use the fact that, under conditions of Proposition 3,
\[
 \frac{b_n}{b_{2n}}\to  \left(\frac{1}{2}\right)^\rho.
 \]
 So, if we know $\lim \frac{b_n}{b_{2n}}$, then we may find $\rho$. In particular, if $ \frac{b_n}{b_{2n}}\to 0$, then $\rho=\infty$, and the distribution $\alpha$ is concentrated at $1$, while  if $ \frac{b_n}{b_{2n}}\to 1$, then $\rho=0$, and the distribution $\alpha$ is concentrated at $0$.

   \item We may deal with a triangular  array, that is, set $c_i=c_{in}$. Then a limiting distribution, if any, may be practically arbitrary. As an example, consider  an integrable, non-negative   function $g(x)$ on $[0,1]$ and set the coefficient  $c_{in}=g(i/n)$. Then, the limiting distribution will be that with the density
    \[
                                 f(x)=\frac{g(x)}{\int_0^1g(x)dx}.
    \]
    A proof of (\ref{alpha}) in this case  may run similarly to what we did above. Note that when considering a counterpart of (\ref{x_n-example}), we may take into account that in this case
    \[
    b_n=:\sum_{i=1}^n c_{in}\sim n\cdot\int_0^1 g(x)dx .
    \]

{ We consider the case of triangular arrays in more detail in a further publication.}

    \item  Clearly, $\arg \max\{X_1,\dots.X_n\}\stackrel{d}{=}\arg \max\{\tilde{X}_1,\dots.\tilde{X}_n\}$, where $\tilde{X}_i=f(X_i)$, and $f(x)$ is a continuous strictly increasing function. It is easy to verify that the corresponding function $\tilde{r}(x)=r(f^{-1}(x))$. This is a way to ``improve''   $r(x)$.
    \item In the case where the distributions of the $X$'s are not continuous, the above technique needs to be improved.  Regarding  the fact that in this case there may be several  ``winners", one can conjecture that the situation may be fixed if we select from winners one   at random (throw lots). On the other hand, in this case  probability $p_{in}\neq 1/n$ even if the $X_i$'s are identically distributed.
       Consider the simplest\\
      \textbf{ Example}.  Let all $X_i=1$ or $0$ with probabilities $p$ and $q$, respectively.
Then
 \[
 p_{in}=p\cdot1+q\cdot q^{n-1}=p+q^n,
 \]
 However, in the case of selecting a winner at random, $ p_{in}=1/n$  just by symmetry, though the same may be also proved directly.
  \end{enumerate}

 We thank professor Vadim Ponomarenko (SDSU) a bygone conversation  with whom
helped us  to come to the statement of this paper  problem. The problem we discussed with professor Ponomarenko,  concerned the above application example in Section \ref{subsecprop}.

\end{document}